\documentclass[12pt]{iopart}
\usepackage{iopams}
\usepackage{setstack}
\usepackage{graphicx}
\usepackage{subfig}
\usepackage{subeqn}
\usepackage{dcolumn}
\usepackage{nicefrac}
\usepackage{citesort}
\usepackage{epstopdf} 
\newcommand{\diff}[2]{\frac{\mathrm{d}#1}{\mathrm{d}#2}}

\begin{document}

\title{Exploring the adaptive voter model dynamics with a mathematical triple jump}

\author{Holly Silk$^1$, G\"uven Demirel$^2$\footnote{Current affiliation: Nottingham University, Business School, Nottingham, UK}, Martin Homer$^1$, and Thilo Gross$^1$}
\address{$^1$ University of Bristol, Department of Engineering Mathematics and Bristol Centre for Complexity Sciences, Bristol, UK}
\address{$^2$ Max-Planck-Institute for the Physics of Complex Systems, N\"othnitzer Str. 38, Dresden, Germany}

\begin{abstract}
Progress in theoretical physics is often made by the investigation of toy models, the model organisms of physics, which provide benchmarks for new methodologies. For complex systems, one such model is the adaptive voter model. Despite its simplicity, the model is hard to analyse. Only inaccurate results are obtained from well-established approximation schemes that work well on closely-related models. We use this model to illustrate a new approach that combines a) the use of a heterogeneous moment expansion to approximate the network model by an infinite system of ordinary differential equations, b) generating functions to map the ordinary differential equation system to a two-dimensional partial differential equation, and c) solution of this partial differential equation by the tools of PDE-theory. Beyond the adaptive voter models, the proposed approach establishes a connection between network science and the theory of partial differential equations and is widely applicable to the dynamics of networks with discrete node-states.      
\end{abstract}

\maketitle

\section{Introduction}
A core challenge in statistical physics is to understand emergence far from equilibrium. In this area much recent progress 
has been made with models that describe a given complex system as a network \cite{networksbook,DorogovtsevExplosivePercolation,StanleyMultilayerNets}. A network model reduces the system to a set of discrete nodes connected by discrete links. It thus simplifies the constituents of a system but explicitly captures the complex pattern of their interactions from which many system-level properties emerge. Networks thereby provide simplified models of the systems, without simplifying the complexity away. 

Because network models still contain the complex topology of the system, this complexity has to be dealt with in the model analysis. While network models can be explored efficiently numerically \cite{ZschalerBioinformatics}, an appealing feature is that analytical progress can be, and has been, made \cite{networksbook,BarabrasiControl,StanleyMultilayerNets,DoPRL}. A particular class of models where much recent progress has been made are so-called adaptive or coevolutionary networks \cite{adaptiverev,gross2009adaptive}. In these system internal states of the network nodes are subject to dynamics on the network, while the network topology itself evolves, depending on the node states. Thus dynamics on and of the network form a feedback loop that gives rise to various forms of adaptation and self-organization \cite{adaptiverev}.          

Two simple toy systems are provided by the adaptive SIS model \cite{SIS} and the adaptive voter model \cite{vazquez}. The adaptive SIS model is an extension of the classical SIS model of epidemic diseases \cite{AndersonMayBook}, in which agents that are susceptible to some disease try to avoid infection by cutting or rewiring links to susceptible individuals. Recent analysis of this model and its variants have shed light on the interplay between epidemic processes and population structure \cite{MillerPRE,FunkPNAS,ShawPRE,FeffermanPRE,Zhang2014} but also served as a testbed for new mathematical and numerical methods \cite{marceau,ZschalerBioinformatics}. The adaptive voter model is an extension of the seminal voter model, where agents try to minimise dissent by cutting or rewiring connections to nodes that have other internal states. In recent papers \cite{Holme,vazquez,kimura,motif,motif2,swarming} this model has been studied with a wide variety of techniques \cite{momentrev}, and can thus be considered as a benchmark system.

Although the adaptive voter and the adaptive SIS model are very similar, the adaptive SIS model seems to be an ``easy'' model 
where even simple analytical approaches perform very well. By contrast the adaptive voter model is a ``difficult'' system, where even the most sophisticated current approaches perform relatively poorly in some regions of parameter space \cite{momentrev}.     

In the mathematical exploration of a network model, the central challenge is often to find an approach that maps the network problem onto a tractable set of equations. 
One prominent class of methods for approximating network dynamics are moment expansions \cite{Keeling,Bauch,SIS,Rogers,kimura,momentrev}. 
The central idea is to write evolution equations for the abundance of certain motifs in the network.
Different expansions can be distinguished by the basis of motifs that they use. 
The most basic approximations, the homogeneous mean-field approximations, only track the abundance of motifs consisting of single nodes.
More sophisticated approximations, such as the homogeneous pair approximation \cite{Keeling,Bauch,SIS} or homogeneous triplet approximation \cite{kimura,momentrev} track also the abundance of larger motifs, such as linked pairs or triplets of nodes.
A powerful class of approximations that have been recently proposed \cite{vazquez,pugliese,marceau,gleeson} are the \emph{heterogeneous approximations}. Here, the definition of a certain motif prescribes not only the state and connectivity of the network nodes in the motif, but also the links connecting them to the rest of the network. 

For illustration, consider an epidemic model describing a network of agents which are either infected with, or susceptible to, some disease. In this case a homogeneous mean field equation keeps track of the numbers of infected and susceptible agents, a homogeneous pair approximation additionally keeps track of the number of links between infected agents, the number of links between susceptible agents, and the number of links connecting susceptible to infected agents. By contrast the heterogeneous approximation keeps track of the number of agents with a given state and a given number of neighbours in specific states, for example, the number of susceptible agents who have exactly three susceptible neighbours and exactly two infected neighbours.

Heterogeneous approximations have been shown to yield excellent results in examples \cite{marceau,gleeson}. However, they typically lead to infinite-dimensional systems of ordinary differential equations (ODEs). 
In practice, the number of equations is limited by a truncation, but the number of ODEs retained is still often of the order of $10^6$ \cite{momentrev}. These high-dimensional ODE systems are then typically studied by numerical integration.  

In this paper we use \emph{generating functions} \cite{genfunctions} to map the infinite-dimensional ODE system from a heterogeneous active neighbourhood approximation \cite{marceau} to a partial differential equation (PDE).
This mapping is exact and reversible and thus does not involve additional assumptions. Generating functions are a major tool of discrete mathematics and as such have been applied to network problems. However, they are typically used to capture the structure rather than the dynamics of networks \cite{networksbook}. For the present context, relevant work includes \cite{wieland} where generating functions were used to explore specific processes in the dynamics of an adaptive epidemic model.

In summary, the approach proposed here is reminiscent of a triple jump, where an athlete crosses a distance by a sequence of three jumps, which employs a different technique.  
Using a \emph{heterogeneous moment expansion} we approximate an agent-based model by a high-dimensional system of ODEs. Then we convert the high-dimensional ODEs into a low-dimensional system of partial differential equations PDEs using \emph{generating functions}. Finally, we solve the PDE system, using methods from the literature.
The key outcome is that the resulting PDE can be solved with much less computational effort than any of the original high-dimensional ODE systems, and may admit analytic solutions.      

We illustrate the proposed approach by applying it to the adaptive voter model in two different ways. In sections~\ref{secVoter}-\ref{sec:characteristics} we exploit a symmetry that holds to good approximation in the voter model, but does not generally exist in other adaptive networks. This simplifies the PDE for the generating function and enables a concise presentation of the triple jump methodology. We are able to solve the PDE analytically by the \emph{method of characteristics} \cite{ockendon}, yielding results in good agreement with agent-based simulations. In section.~\ref{sec:fullmodel} we avoid exploiting the symmetry, which results in a somewhat more involved but more general treatment.

Our main intention is to show that the systems of equations obtained from the heterogeneous approximation can be solved relatively straightforwardly using generating functions followed by a suitable PDE technique. The triple jump approach can reveal full time-dependent solutions that describe the network dynamics to a high degree of accuracy. We believe that this approach will be valuable for the wide variety of models in which the heterogeneous approximation affords an almost exact description of the system.   

\section{Adaptive voter model \label{secVoter}}
The adaptive voter model consists of a network of $N$ nodes, representing agents, and $K$ bidirectional links, representing social contacts.
Each agent $i$, is associated with a binary variable $s_i \in \{ \mathrm{A},\mathrm{B} \}$ representing the agent's opinion.
We initialize the network of agents as an Erd\H{o}s-R\'enyi random graph and assign opinions to the agents randomly with equal probability.
The network is then evolved in time by consecutive update steps. 
In each step a link connecting two agents $i$ and $j$ is selected randomly.
If the agents hold identical opinions, $s_i=s_j$, then the link is said to be \emph{inert} and nothing happens.
If $s_i \neq s_j$, then the link is said to be \emph{active} and one agent, $a\in\lbrace i,j\rbrace$, is chosen randomly (with probability 1/2) to resolve the conflict.
With probability $p$, the link connecting $i$ and $j$ is cut and agent $a$ establishes a new link to a randomly chosen agent, $k$, with $s_k=s_a$ (chosen uniformly across such agents).
Otherwise, with probability $\bar{p}=1-p$, agent $a$ changes its opinion such that $s_i=s_j$. 
We say that in the former case the conflict is resolved by a \emph{rewiring event}, and in the latter case by an \emph{opinion adoption event}.

One can verify that the rules of the adaptive voter model are unbiased \cite{momentrev}, such that there is no net drift in the number of nodes holding a particular opinion. Thus, the fraction of nodes holding opinion A remains constant in time, except for stochastic fluctuations that become negligible in the thermodynamic limit $N\to\infty$ (with constant $K/N$).  
For simplicity, we can thus begin by focusing on the symmetric case where the number of nodes holding opinions A and B are equal.

\section{Heterogeneous moment expansion}
\begin{figure}
\centerline{\includegraphics[width=0.5\textwidth]{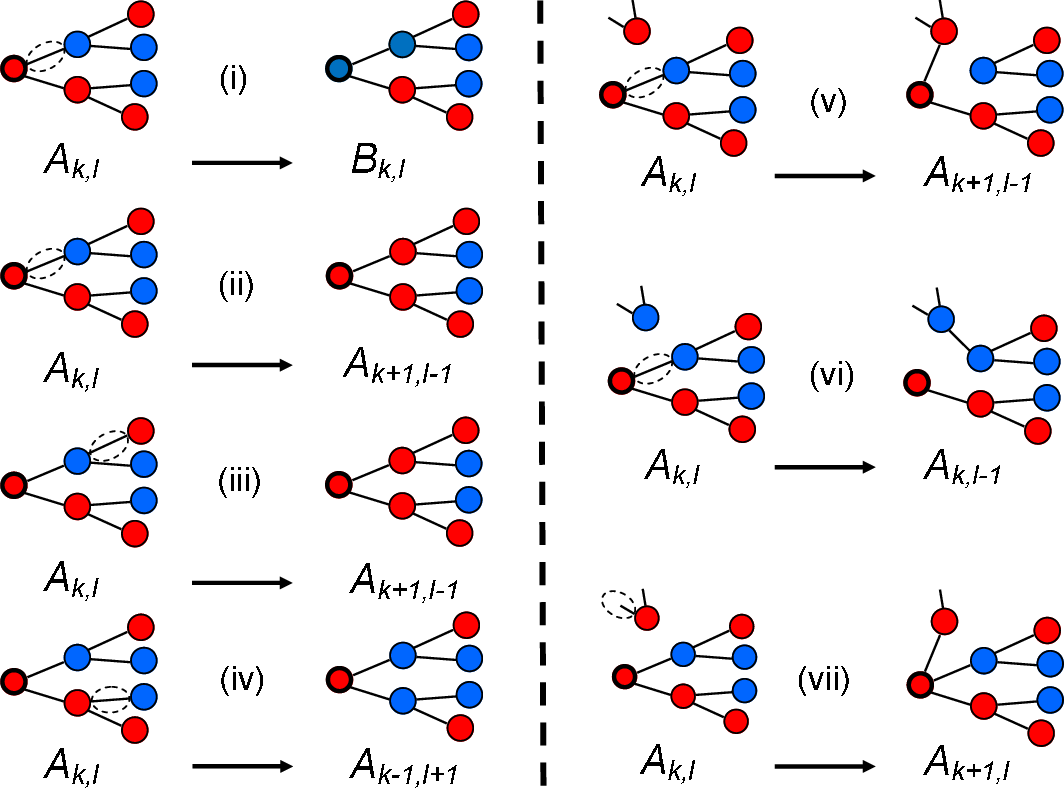}}
\caption{The different update events of the active neighbourhood approximation. Equation (\ref{eq:ODE}) captures the change in the abundance of the focal  node (bold). Colour is used to indicate the state of the respective nodes. The left panel shows opinion adoption events and the right panel rewiring events. The numerical labelling match that of equation (\ref{eq:ODE}).}
\label{fig:link_update}
\end{figure}
In the first step of the triple jump, we convert the stochastic agent-based model into an infinite-dimensional system of ordinary differential equations (ODEs) by a heterogeneous expansion, known as the active neighbourhood approximation.
Following \cite{marceau} we define $A_{k,l}$ (respectively $B_{k,l}$) to be the normalized number of agents of opinion A (respectively B) who have $k$ A-neighbours and $l$ B-neighbours.  In the thermodynamic limit we can treat the $A_{k,l}$ as continuous variables and capture their dynamics by the differential equations
\renewcommand{\thesubequation}{\themainequation\roman{equation}}
\begin{subequations}\label{eq:ODE}
\begin{eqnarray}
	\frac{\rmd A_{k,l}}{\rmd t}
	&= \frac{\bar{p}}{2}\left[ kB_{k,l} - lA_{k,l}\right]  
	\label{eq:Akl-i}\\
	&\quad+ \frac{\bar{p}}{2}\left[ (l+1)A_{k-1,l+1} - lA_{k,l}\right]  
	\label{eq:Akl-ii}\\
 	&\quad+  \frac{\bar{p}}{2}\frac{\sum_{k,l}(k-1)kB_{k,l}}{\sum_{k,l} k 
 	B_{k,l}}\left[ (l+1)A_{k-1,l+1} - lA_{k,l}\right] 
 	\label{eq:Akl-iii}\\
 	&\quad+  \frac{\bar{p}}{2}\frac{\sum_{k,l}lkA_{k,l}}{\sum_{k,l}k 
 	A_{k,l}}\left[ (k+1)A_{k+1,l-1} - kA_{k,l}\right] 
 	\label{eq:Akl-iv}\\
 	&\quad+  \frac{p}{2}\left[ (l+1)A_{k-1,l+1} - lA_{k,l}\right] 
 	\label{eq:Akl-v}\\
 	&\quad+  \frac{p}{2}\left[ (l+1)A_{k,l+1} - lA_{k,l}\right]
 	\label{eq:Akl-vi}\\
 	&\quad+  \frac{p}{2}\frac{\sum_{k,l}lA_{k,l}}{\sum_{k,l}A_{k,l}}\left[ A_{k-1,l} - A_{k,l}\right].
 	\label{eq:Akl-vii}
\end{eqnarray}
\end{subequations}
The terms in equation \eref{eq:ODE} describe the change experienced by a focal node $A_{k,l}$ due to opinion adoption (i to iv) and rewiring (v to vii) events. Specifically, contributions arise from 
(i) the focal node adopting the opinion of a neighbour, 
(ii) a neighbour of type B adopting the opinion of the focal node, 
(iii) a neighbour of type B adopting the opinion of another node of type A, 
(iv) a neighbour of type A adopting the opinion of a node of type B, 
(v) the focal node rewiring one of its links away from a neighbour of type B (acquiring a new neighbour of type A), 
(vi) a neighbour of type B rewiring a link away from the focal node, 
and (vii) a node of type A rewiring one of its links to the focal node. 

Let us consider the first term (i) in more detail. Adopting the opinion of a neighbour can result in both a loss or creation of $A_{k,l}$ nodes. A loss occurs when $A_{k,l}$ nodes are convinced to adopt the opinion of a neighbour of type B (hence becoming $B_{k,l}$ nodes); this occurs at a rate proportional to $A_{k,l}$ and to the amount of B neighbours of $A_{k,l}$, namely $l$. Similarly, $A_{k,l}$ nodes are created when $B_{k,l}$ nodes are convinced by one of their $k$ A neighbours, at a rate proportional to $kB_{k,l}$. The probability of an opinion-adoption event occurring is $\bar{p}/2$.
Thus the contribution of events of type (i) to the ODE is
\begin{math}
  \frac{\bar{p}}{2}\left[ kB_{k,l}-lA_{k,l} \right].
\end{math}

The events (iii), (iv), and (vii) involve nodes outside the direct neighbourhood of the focal node.
They are thus dependent on the number of next-nearest neighbours (iii, iv) or active links existing elsewhere (vii). The corresponding rates then depend on longer-ranged correlations that are not captured by 
the $A_{k,l}$ and $B_{k,l}$ alone. We therefore need to estimate these rates based on the available information from the nearest-neighbour correlations captured. This approximation is called moment closure, and is known to be the main source of inaccuracy in heterogeneous moment expansions for networks \cite{momentrev}.

For instance, the rate (iv) at which a typical neighbour of type A (A-neighbour) of the focal node adopts the opinion B depends on the average number of B-neighbours of the A-neighbour, and thus on the next-nearest-neighbourhood of the focal node. To approximate this number, one considers all potential A-neighbours, based on the known distribution $A_{k,l}$, but takes into account that a node that has $k$ links to other A-nodes is $k$ times more likely to be an A-neighbour of the focal node. The distribution of A-neighbours of the focal node is thus $ k A_{k,l}/(\sum k A_{k,l})$, where the denominator arises from normalization. Based on this distribution we can then estimate the number of B-neighbours of a typical A-neighbour of the focal node as $(\sum lk A_{k,l})/(\sum k A_{k,l})$, which appears as a factor in the corresponding term (iv).

So far the ODE system (1), does not constitute a closed model because the variables $B_{k,l}$ appear in the equations. Symmetry can be exploited in two different ways to deal with these variables. First, we can treat the $B_{k,l}$ as genuine dynamical variables, which follow a set of differential equations that are symmetric to the equations governing $A_{k,l}$. While this approach does not involve any additional assumptions, it gives rise to a two-dimensional PDE, which is not straightforward to solve in closed form. Second, we can exploit symmetry more directly by finding a suitable approximation that eliminates the $B_{k,l}$ from the equation. Here, we use $(\sum k(k-1) B_{k,l})/\sum k B_{k,l}) = (\sum l(l-1) A_{k,l})/\sum l A_{k,l})$, and $lA_{k,l} \approx kB_{k,l}$. The first of these relationships states that the number of $ABA$ triplets per $AB$ link is equal to the number of $BAB$ triplets per $BA$ link. This relationship is exact in the statistical sense and thus does not constitute an additional assumption. The second relationship can be motivated from results of the homogeneous expansion \cite{momentrev} but is more problematic as it differs from the statistically exact $ lA_{k,l} =  l B_{l,k}$. We nevertheless use this approximation because the exact expression would give rise to a non-local PDE which is more difficult to solve than the full two-dimensional PDE system. 

In the following we pursue both of the routes outlined above. First, in sections~\ref{sec:genfun} and \ref{sec:characteristics}, we investigate the simplified system using the approximation $lA_{k,l} \approx kB_{k,l}$, which leads to a one-dimensional PDE that is then solved analytically. Then, in section~\ref{sec:fullmodel}, we return to the full systems where the $B_{k,l}$ are treated as additional variables. This leads to a two-dimensional PDE for which an analytical solution is conceivable, but which we solve here using a highly efficient numerical step.
 
\section{Generating functions}\label{sec:genfun}
Generating functions \cite{genfunctions} can be used to reduce the infinite-dimensional ODE system to a finite-dimensional PDE by interpreting the $A_{k,l}$ as Taylor-like coefficients of a polynomial, whose time evolution is then studied.
We start by defining the generating function $Q(t,x,y)= \sum_{k,l}A_{k,l}(t)x^{k}y^{l}$, where $x$ and $y$ are abstract spatial variables that do not have any physical interpretation, but act as an indexing mechanism. Finding $Q(t,x,y)$ is equivalent to finding all of the moments $A_{k,l}$ and thus constitutes a solution of the system.

The time evolution of $Q$ is given by 
\begin{equation}
\frac{\partial Q}{\partial t} = Q_t = \sum_{k,l} \frac{\rmd A_{k,l}}{\rmd t}  x^{k}y^{l}.
\end{equation} 
Substituting \eref{eq:ODE} we obtain an expression in which the right-hand side can be written again in terms of the function $Q$ and its derivatives.
For instance, the process (ii), described above, results in a term proportional to 
\begin{equation}
\sum_{k,l} \left[ (l+1)A_{k-1,l+1} - lA_{k,l}\right] x^{k}y^{l} = xQ_y - yQ_y,
\end{equation}
where $Q_x=\partial Q / \partial x$. The validity of this identity can be verified by separating the two terms in the square bracket and shifting the indices $k-1 \to k$, $l+1 \to l$ on the first.
Proceeding analogously, we find that $Q$ satisfies
\begin{equation}
  Q_t= \frac{\bar{p}\beta}{2}(y-x)Q_x + \left[ \left(\frac{1+\bar{p}\alpha}{2}\right)(x-y) + \frac{p}{2}(1-y)\right] Q_y + \frac{p\gamma}{2}(x-1)Q,
\label{eq:pde}
\end{equation}
where
\begin{equation}
  \alpha = \frac{Q_{yy}(t,1,1)}{Q_{y}(t,1,1)},\qquad
  \beta = \frac{Q_{xy}(t,1,1)}{Q_{x}(t,1,1)}, \qquad
  \gamma = \frac{Q_{y}(t,1,1)}{Q(t,1,1)}
\label{eq:alphabetagamma}
\end{equation}
are the transformed factors from the moment closure approximation. They represent the density of $ABA$ triplets per $AB$ link, the density of $AAB$ triplets per $AA$ link and the density of $AB$ links per $A$ nodes, respectively. For the moment we will treat these factors as unknown parameters, but determine them later from a self-consistency condition.

\section{Solving the PDE}\label{sec:characteristics}
The generating function PDE \eref{eq:pde} is a first-order scalar quasilinear equation, of the form
\begin{equation}
  a(t,x,y)Q_t+b(t,x,y)Q_x+c(t,x,y)Q_y=d(t,x,y,Q).
\label{eq:pde-general}
\end{equation}
Such PDEs can be solved by the method of characteristics \cite{ockendon}. The central 
idea is to describe the solution surface $Q=Q(t,x,y)$ parametrically.
Three parameters are required, labelled $\eta$, $\xi_1$, $\xi_2$.
The method captures the dynamics with a low-dimensional set of ODEs:
the \emph{bicharacteristic equations}
\begin{equation*}
  \diff{t}{\eta}=a,\quad \diff{x}{\eta}=b,\quad \diff{y}{\eta}=c,\quad \diff{Q}{\eta}=d.
\label{eq:bichars-general}
\end{equation*}
 The family of solutions of the bicharacteristic equations, indexed by parameters
$(\xi_1,\xi_2)$, makes up the solution surface.

One can conceptualise this construction, by writing \eref{eq:pde-general} in the vector form
\begin{equation*}
  (a, \, b, \, c,\, d)\cdot\left(Q_t, \,
    Q_x, \, Q_y, \, -1 \right)=0.
\end{equation*}
Since the second vector is normal to the solution surface  $Q=Q(t,x,y)$,  in $(t,x,y,Q)$-space, the vector $(a,b,c,d)$ is everywhere tangent to the solution surface. Hence curves $(t(\eta),x(\eta),y(\eta),Q(\eta))$ that satisfy $(\rmd t/\rmd\eta,\rmd x/\rmd\eta,\rmd y/\rmd\eta,\rmd Q/\rmd\eta)=(a,b,c,d)$ remain on the solution surface for all $\eta$.

For the specific system (\ref{eq:pde}) the bicharacteristic equations are
\begin{eqnarray}
\frac{\rmd t}{\rmd \eta} &= 1, \label{eq:char-t} \\
\frac{\rmd x}{\rmd \eta} &= \frac{\bar{p}\beta}{2}(x-y), \label{eq:char-x}\\
\frac{\rmd y}{\rmd \eta} &= \frac{1}{2}(1+\bar{p}\alpha)(y-x) + \frac{p}{2}(y-1), \label{eq:char-y}\\
\frac{\rmd Q}{\rmd \eta} &= \frac{p\gamma}{2}(x-1)Q \label{eq:char-Q}.
\end{eqnarray}
The quantities $\alpha$, $\beta$, $\gamma$ can be regarded as auxiliary variables that change dynamically in time. However, as we are primarily interested in the long term behaviour, where also these variables become stationary, we can treat them as parameters and determine their steady state values later by demanding self-consistency.  

To express the solutions of \eref{eq:char-t}--\eref{eq:char-Q} it is useful to define the matrix $M$ of the homogeneous linear operator defined by \eref{eq:char-x} and \eref{eq:char-y};
\begin{equation}
M = \frac{1}{2}\left( \begin{array}{cc} \bar{p}\beta & - \bar{p}\beta \\
-(1+\bar{p}\alpha) & 1 + p + \bar{p}\alpha \end{array}\right) \label{eq:M}.
\end{equation}  
In the following, we denote the eigenvalues of this matrix as $\lambda_{1}$,$\lambda_{2}$ and the corresponding eigenvectors as $\left[ v_{1}^{1}\; v_{1}^{2}\right]^{\rm T}$, $\left[ v_{2}^{1}\; v_{2}^{2}\right]^{\rm T}$.

Since \eref{eq:char-x} and \eref{eq:char-y} are independent of $t$ and $Q$, solving for $x$ and $y$ is straightforward,
giving expressions in terms of the eigenvalues and eigenvectors of \eref{eq:M}.
Then, solving equations \eref{eq:char-t} and \eref{eq:char-Q} results in an analytic solution for $Q$
\begin{eqnarray}\label{eq:Qsoln}
\fl	Q = \exp \left\lbrace\frac{\langle k\rangle}{2}\left( \xi_1 + \xi_2 -2\right) - \frac{p\gamma}{2 \left(v^{1}_{1}v_{2}^{2}- v_{1}^{2} v_{2}^{1}\right)} \left[ \frac{v^{1}_{1}\left( v_{2}^{2}\left( \xi_1 - 1 \right) - v_{2}^{1} \left(\xi_2 - 1 \right)\right)}{\lambda_1}\left(\rme^{\lambda_{1}t}-1\right) \right.\right.\nonumber \\
\lo+ \left.\left.\frac{v_{2}^{1} \left( v_{1}^{1} \left( \xi_2 - 1 \right) - v_{1}^{2}\left( \xi_1 - 1 \right)\right)}{\lambda_2}\left(\rme^{\lambda_{2}t}-1\right)\right]\right\rbrace,
\end{eqnarray}
where $(\xi_1,\xi_2)$ parametrise the initial state of the network. The solution \eref{eq:Qsoln} can then be written in terms of $x$ and $y$, using the solutions to \eref{eq:char-x} and \eref{eq:char-y}.

From this closed form solution for the generating function $Q$ all properties of the distribution $A_{k,l}$ can be computed analytically.
Moreover, we can directly compute aggregate properties such as the number of active links $\sum_{k,l} lA_{k,l} = Q_{y}(1,1)$.
In particular we can now obtain consistency conditions for $\alpha$, $\beta$, and $\gamma$, using their definition \eref{eq:alphabetagamma}.
Combining the individual conditions we find $\alpha=\beta=\gamma$.
Above we defined $\gamma$, as the number of active links per node of type A, $\alpha$ is the number of $ABA$ triplets per $AB$-link, and $\beta$ is the number of $AAB$ triplets per $AA$ link. In other words, the three quantities can be interpreted as the expected number of active links found by picking a random node of type $A$ ($\gamma$), following a random $AB$-link ($\alpha$), and following a random $AA$-link ($\beta$). The identity of these three expectations thus implies the absence of correlations beyond the nearest neighbour interactions. This is clearly an artefact of our assumptions, as it is known that longer ranged correlations play a role in the adaptive voter model \cite{momentrev}.

\begin{figure}
\centerline{\includegraphics[width=0.75\textwidth]{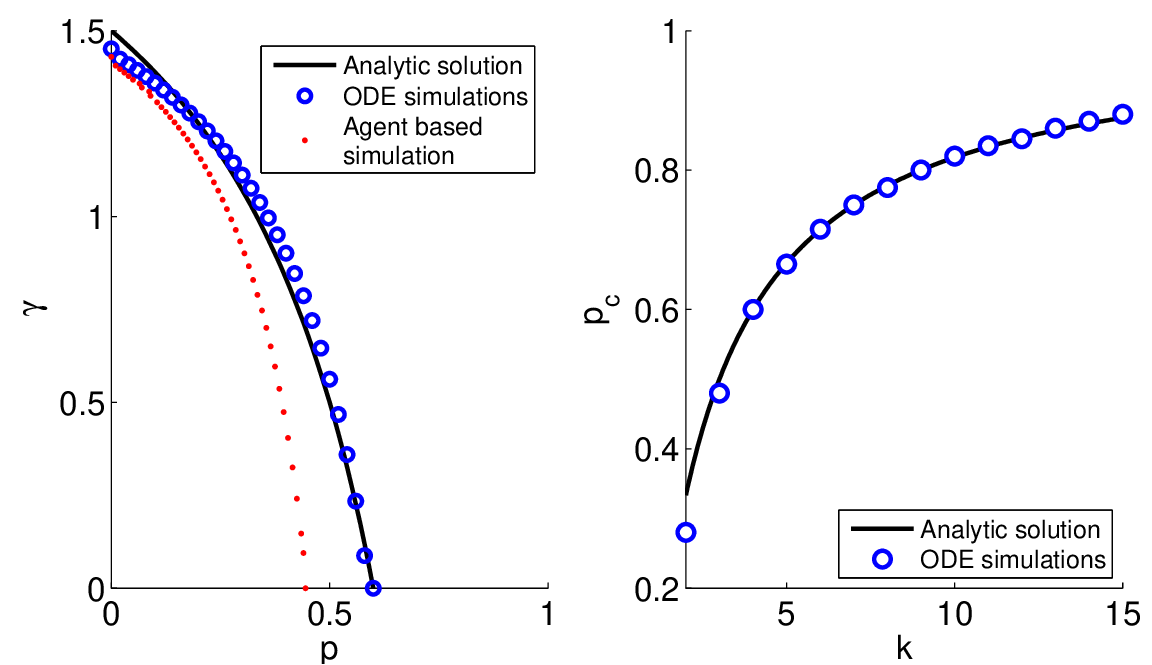}}
\caption{Fragmentation transition. Shown (left (a)) are number of active links $\gamma(\infty)$ as a function of the rewiring rate $p$ (Parameters: $\langle k \rangle =4$, for agent based simulation $N=10^5$, for ODE simulation $k_{\rm max}=1000$ except at transition point where $k_{\rm max}=100$), and (right (b)) critical rewiring rate $p_{\rm c}$ as a function of mean degree $\langle k\rangle$.}
\label{fig:gammavsp}
\end{figure}

The analytical results recapture the well-known behaviour of the adaptive voter model (figure~ \ref{fig:gammavsp}).
If the rewiring rate $p$ is sufficiently low then the system approaches an active state ($\gamma\neq0$) that is stable in the thermodynamic limit \cite{vazquez}. The system remains in this active state for a long time (the value $\gamma$ in figure~\ref{fig:gammavsp}(a)) before achieving consensus in one opinion \cite{Rogers2013}. If the rewiring rate exceeds a threshold $p_{\rm c}$, however,  the system rapidly approaches a fragmented state ($\gamma=0$), where the network breaks into two disconnected components that hold different opinions, but are internally in consensus. The figure shows that the analytical solution is in good agreement with results from the numerical integration of the high dimensional ODE system from the heterogeneous expansion.

We can also analytically determine the critical rewiring rate, $p_{\rm c}$. Using the result $\alpha = \beta = \gamma$ we find that the values of the eigenvalues and eigenvectors of \eref{eq:M} are $\lambda_{1,2} = \left( 1+p \right)/4 + \bar{p}\gamma/2 \pm \sqrt{d}/2$, $v_1^1,v_2^1=1$ and $v_1^2,v_2^2 = \left[ \left( 1+p \right)/2 \mp \sqrt{d}\right] /(\bar{p}\gamma)$
where $d = \left( 1+p\right) ^2/4 + \gamma^2\bar{p}^2 + \gamma\bar{p}$.
The eigenvalues $\lambda_1,\lambda_2$ are positive except at the special cases of $p=1$, $p=0$, where one eigenvalue is equal to zero. Therefore, ignoring the special cases, any exponential terms in $Q_x(1,1)$ approach zero as $t\rightarrow\infty$ and we are left with a constant, $Q_x\left( 1,1\right)= \gamma + (1+p)/(1-p)$. By definition $Q_y(1,1) = \gamma$. Substituting these results into $Q_x(1,1)+ Q_y(1,1) = \langle k\rangle$ gives
$\alpha = \beta = \gamma = \left[ \langle k\rangle - (1+p)/(1-p)\right]/2$.
At the fragmentation transition the active links vanish, thus $\gamma=0$ and $p_{\rm c} = (\langle k\rangle - 1)/(\langle k\rangle + 1)$.
Because the longer-ranged correlations are not captured this result is only qualitatively correct, as shown in figure~\ref{fig:gammavsp}(a). Displayed are the results from agent-based simulations, where the number of agents $N = 10^5$, as well as those of the analytic solution \eref{eq:Qsoln} and numerical integration of the high-dimensional ODE system \eref{eq:ODE}. There is good agreement for small values of $p$, however the moment expansion (in both cases) fails to quantitatively capture the behaviour of the model close to the fragmentation point, overestimating $p_{\rm c}$. This is in line with other moment expansion techniques which also tend to overestimate the fragmentation point \cite{momentrev}. The active neighbourhood approximation nevertheless performs better than other simpler approximations \cite{momentrev}, with the exception of \cite{motif} which is only applicable close to the fragmentation point.

\section{Full model}\label{sec:fullmodel}
Previous work and also our results from the previous section show that the active neighbourhood approximation does not yield good results for the adaptive voter model at low mean degree. The approach proposed here is nevertheless useful at high mean degree and, more importantly, also for the large class of other similar adaptive network models that have been proposed. In our treatment above, we have exploited a specific symmetry of the adaptive voter model by assuming $lA_{k,l} = kB_{k,l}$ and cancelling the terms of (\ref{eq:Akl-i}). This simplification is not generally possible for other models, such as the adaptive SIS model. We therefore now repeat our analysis of the voter model, without the simplifying assumptions.

We begin by using the active neighbourhood approximation to write rate equations for the two variables $A_{k,l}$, $B_{k,l}$ which produces two infinite-dimensional systems of ODEs. The rate equation for $A_{k,l}$ is given by (\ref{eq:ODE}) with a similar equation found for $B_{k,l}$. Then by defining the generating functions $Q(t,x,y) = \sum_{k,l}A_{k,l}(t)x^ky^l$ and $R(t,x,y) = \sum_{k,l}B_{k,l}(t)x^ky^l$ we arrive at the following pair of coupled PDEs
\begin{eqnarray}\label{eq:coupled_PDES}
\fl \left( \begin{array}{c}Q_t\\
 R_t
 \end{array}\right) 
 =\left( \begin{array}{cc}\frac{\bar{p}\beta}{2}(y-x) & \frac{\bar{p}}{2}x \\ 0 & \left[\frac{\bar{p}}{2}(\delta+1)+ \frac{p}{2}\right](y-x)-\frac{x}{2} + \frac{p}{2}
\end{array}\right) 
\left( \begin{array}{c}Q_x\\
R_x
\end{array}\right)\nonumber
\\
+ \left( \begin{array}{cc}\left[\frac{\bar{p}}{2}(\alpha+1)+ \frac{p}{2}\right](x-y)- \frac{y}{2} + \frac{p}{2}&0 \\ \frac{\bar{p}}{2}y & \frac{\bar{p}\varepsilon}{2} (x-y)
\end{array}\right) 
\left( \begin{array}{c}Q_y\\
R_y
\end{array}\right)\nonumber
\\
+ \left(\begin{array}{cc}\frac{p\gamma}{2}(x-1) & 0 \\
0 & \frac{p\zeta}{2}(y-1)
\end{array}\right)
\left(\begin{array}{c}Q\\
R
\end{array}\right),
\end{eqnarray}
where
\begin{eqnarray}\label{eq:coupled_consistency}
\alpha &=\frac{R_{xx}(t,1,1)}{R_x(t,1,1)},\qquad \beta &=\frac{Q_{xy}(t,1,1)}{Q_x(t,1,1)},\qquad \gamma =\frac{Q_y(t,1,1)}{Q(t,1,1)}\nonumber\\
\delta &=\frac{Q_{yy}(t,1,1)}{Q_y(t,1,1)},\qquad \varepsilon &=\frac{R_{xy}(t,1,1)}{R_y(t,1,1)},\qquad \zeta =\frac{R_{x}(t,1,1)}{R(t,1,1)}.
\end{eqnarray}

We consider the steady state where $\partial/\partial t = 0$ and, to solve for the generating functions $Q$ and $R$, expand them as Taylor series, since the extension of the method of characteristics to PDE systems of this nature is non-trivial. We expand $Q$ and $R$ around the point $x=y=1$. By substituting these expressions into \eref{eq:coupled_PDES} we can equate coefficients of powers of $x$ and $y$ and solve a set of simultaneous equations for the derivatives of $Q$, $R$ at $x = y = 1$. Then the consistency conditions (\ref{eq:coupled_consistency}) along with the conservation of links $Q_x(1,1)+Q_y(1,1)+R_x(1,1)+R_y(1,1)=  \langle k\rangle$, and nodes $Q(1,1)+R(1,1)=1$,  can be used to find the abundance of active links in the system ($\gamma$, $\delta$) in the steady state along with other motifs (\ref{eq:coupled_consistency}). This can, in principle, be used to find all terms in the Taylor series for $Q$ and $R$ and, as we show below, good results are obtained already by considering the first 7 terms.

\begin{figure}
\centerline{\includegraphics[width=0.75\textwidth]{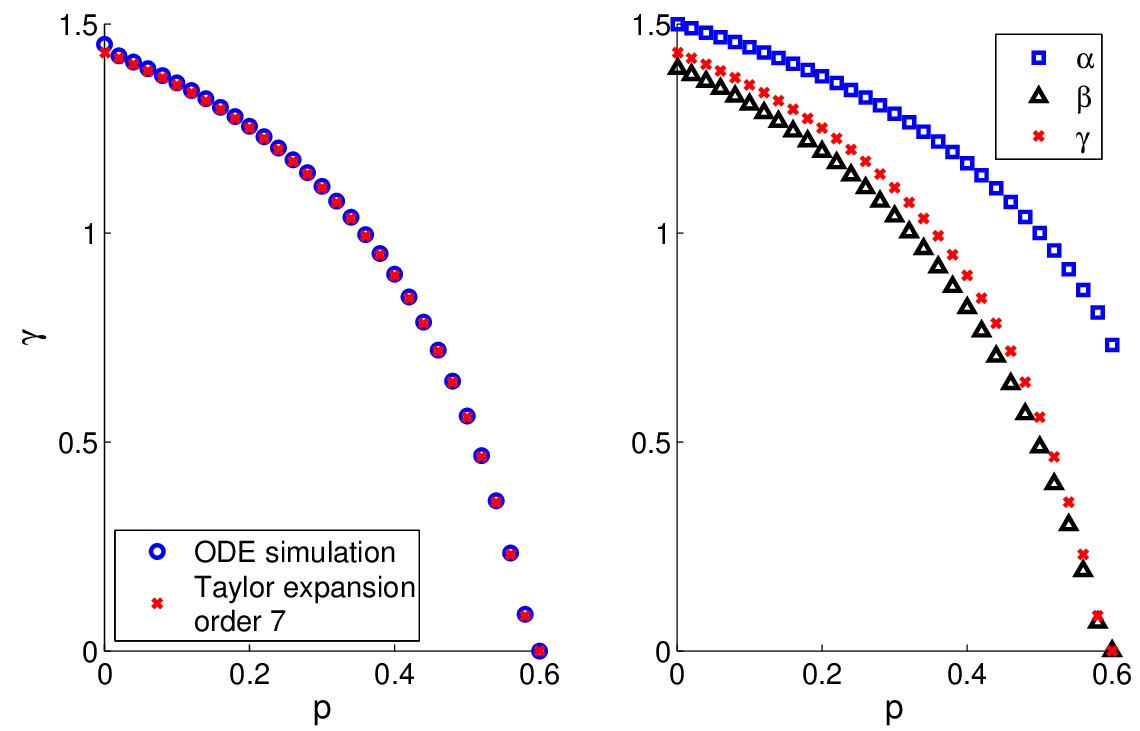}}
		\caption{Fragmentation transition. Shown are number of active links $\gamma(\infty)$ as a function of the rewiring rate $p$ (left (a), for ODE simulation $k_{\rm max}=1000$ except at transition point where $k_{\rm max}=100$), where the results of the Taylor expansion are in good agreement with the ODE simulation. Right (b) are the values of $\alpha$, $\beta$ and $\gamma$ obtained from the Taylor expansion. Mean degree $\langle k \rangle = 4$.}
    	\label{fig:taylor_results}
\end{figure}

Figure~\ref{fig:taylor_results}(a) shows the results obtained from the seventh-order Taylor expansion of $Q$ and $R$ in the symmetric case, where the densities of $A$s and $B$s in the system are equal, compared with those from numerical simulation of the high-dimensional ODEs \eref{eq:ODE}. For this symmetric case we find that $\alpha= \delta$, $\beta = \varepsilon$ and $\gamma = \zeta$. The Taylor expansion results are in excellent agreement with the ODE simulations for all values of $p$, and an improvement on the analytic results in the previous section, illustrating the information that is lost by cancelling the terms (\ref{eq:Akl-i}), even in the symmetric case. This is also apparent in figure~\ref{fig:taylor_results}(b) which shows that in the full model solved here, $\alpha \neq\beta \neq \gamma$. We also note the good agreement between our results from the previous section and the ODE simulation, despite a significant error in the value of $\alpha$.
Furthermore, figure~\ref{fig:taylor_results}(a) demonstrates that very little information is lost from the heterogeneous expansion by truncating the Taylor series solution for $Q$ and $R$ to a relatively small number of terms. It is therefore possible to obtain excellent results without having to solve the full ODE system; we have produced the same results as the numerical integration of the high-dimensional set of ODEs  \eref{eq:ODE}  of the heterogeneous expansion, while avoiding the computational costs associated with the numerical integration \cite{momentrev}.

\section{Conclusions}
In this paper we have proposed an approach for the investigation of network dynamics that combines heterogeneous expansions, generating functions, and solution of the resulting PDE. 
Using this mathematical triple jump, analytic solutions for heterogeneous expansions can be obtained, which 
we demonstrated on the example of the adaptive voter model. 
The triple jump approach generally does not involve critical assumptions, other than those inherent in the heterogeneous expansion, as demonstrated by our second method. By looking at the whole system we have shown that it is possible to obtain accurate results from a heterogeneous moment expansion without having to simulate the full system of ODEs. 

Here, we focused on the adaptive voter model because its symmetry allows for a concise presentation of the triple jump methodology. When applying the approach to the full model we chose to represent states as a variable ($B_{k,l}$ in addition to $A_{k,l}$). Another possibility would be to instead have the state as an extra index ($A_{k,l,m}$ instead of $A_{k,l}$ where the index $m$ signals the state of the focal node.), leading to a single scalar PDE in a three-dimensional space, rather than two coupled PDEs in a two-dimensional space. The relative benefits of each method will vary depending on the resulting PDEs and are thus model dependent.

We observed that the additional assumption $lA_{k,l}=kB_{k,l}$ which we made in the simplified treatment of the system, has the  unexpected effect of removing the link correlations from the active steady state. This is interesting because similar correlations are the main cause of inaccuracies in moment expansions. A detailed investigation of why the additional assumption causes these correlations to vanish could yield new insights into the emergence of correlations in general adaptive network models. Thereby it could lead to the discovery of approximation schemes that accommodate such correlations better than current approaches.

For the adaptive voter model, it is known that the heterogeneous approximation we utilise here provides only qualitative results \cite{momentrev}. However, heterogeneous expansions provide an excellent approximation in other models \cite{marceau,gleeson}. We hope our methodology will be a useful tool in these instances. Beyond analytical solution of the PDEs for simple systems, the perspective that is opened up here is to transfer insights from PDE theory to the analysis of adaptive networks. Such insights concern for instance results on information flow and uniqueness of solutions and may thus lead to a deeper understanding of network dynamics.

\ack
This work was partly funded by the EPSRC through the Bristol Centre for Complexity Sciences (BCCS) and grant EP/K031686/1.

\section*{References}
\bibliographystyle{iopart-num}
\bibliography{triple_jump_references}

\end{document}